\documentclass[a4paper,12pt]{amsart}
\usepackage{amsfonts}
\usepackage{amssymb}
\usepackage{ifthen}
\usepackage{amscd}
\usepackage{amsxtra}
\usepackage{graphicx}
\usepackage{color}
\nonstopmode \numberwithin{equation}{section}
\newtheorem{thm}{Theorem}
\newtheorem{lem}{Lemma}
\newtheorem{cor}{Corollary}

\newtheorem{cl}{Claim}
\newtheorem{ca}{Case}
\newtheorem{sca}{Subcase}
\newtheorem{scl}{Subclaim}
\newtheorem{conj}{Conjecture}

\theoremstyle{definition}
\newtheorem{defn}{Definition}

\newtheorem{op}[equation]{Open Problem}
\newtheorem{ques}[equation]{Question}
\newtheorem{rem}{Remark}[section]
\newtheorem{exam}[equation]{Example}

\newcounter {own}
\def\theown {\thesection       .\arabic{own}}

\newenvironment{pf}[1][]{%
 \vskip 3mm
 \noindent
 \ifthenelse{\equal{#1}{}}%
  {{\slshape Proof. }}%
  {{\slshape #1.} }%
 }%
{\qed\bigskip}

\newcounter{alphabet}
\newcounter{tmp}
\newenvironment{Thm}[1][]{\refstepcounter{alphabet}%
\bigskip%
\noindent%
{\bf Theorem \Alph{alphabet}}%
\ifthenelse{\equal{#1}{}}{}{ (#1)}%
{\bf .} \itshape}{\vskip 8pt}

\makeatletter
\newcommand{\Ref}[1]{\@ifundefined{r@#1}{}{\setcounter{tmp}{\ref{#1}}\Alph{tmp}}}
\makeatother

\newenvironment{Lem}[1][]{\refstepcounter{alphabet}%
\bigskip%
\noindent%
{\bf Lemma \Alph{alphabet}}%
{\bf .} \itshape}{\vskip 8pt}

\newcommand{\ID}{{\mathbb D}}

\def\be{\begin{equation}}
\def\ee{\end{equation}}

\newcommand{\bee}{\begin{enumerate}}
\newcommand{\eee}{\end{enumerate}}

\newcommand{\blem}{\begin{lem}}
\newcommand{\elem}{\end{lem}}
\newcommand{\bthm}{\begin{thm}}
\newcommand{\ethm}{\end{thm}}
\newcommand{\bcor}{\begin{cor}}
\newcommand{\ecor}{\end{cor}}
\newcommand{\beg}{\begin{exam}}
\newcommand{\eeg}{\end{exam}}
\newcommand{\begs}{\begin{examples}}
\newcommand{\eegs}{\end{examples}}
\newcommand{\bdefe}{\begin{defn}}
\newcommand{\edefe}{\end{defn}}
\newcommand{\bprob}{\begin{prob}}
\newcommand{\eprob}{\end{prob}}
\newcommand{\bques}{\begin{ques}}
\newcommand{\eques}{\end{ques}}
\newcommand{\bei}{\begin{itemize}}
\newcommand{\eei}{\end{itemize}}
\newcommand{\bcon}{\begin{conj}}
\newcommand{\econ}{\end{conj}}
\newcommand{\bop}{\begin{op}}
\newcommand{\eop}{\end{op}}

\newcommand{\bca}{\begin{ca}}
\newcommand{\eca}{\end{ca}}
\newcommand{\bsca}{\begin{sca}}
\newcommand{\esca}{\end{sca}}

\newcommand{\bcl}{\begin{cl}}
\newcommand{\ecl}{\end{cl}}

\newcommand{\bscl}{\begin{scl}}
\newcommand{\escl}{\end{scl}}

\newcommand{\bcons}{\begin{conjs}}
\newcommand{\econs}{\end{conjs}}
\newcommand{\bprop}{\begin{propo}}
\newcommand{\eprop}{\end{propo}}
\newcommand{\br}{\begin{rem}}
\newcommand{\er}{\end{rem}}
\newcommand{\brs}{\begin{rems}}
\newcommand{\ers}{\end{rems}}
\newcommand{\bo}{\begin{obser}}
\newcommand{\eo}{\end{obser}}
\newcommand{\bos}{\begin{obsers}}
\newcommand{\eos}{\end{obsers}}
\newcommand{\bpf}{\begin{pf}}
\newcommand{\epf}{\end{pf}}
\newcommand{\ba}{\begin{array}}
\newcommand{\ea}{\end{array}}
\newcommand{\beq}{\begin{eqnarray}}
\newcommand{\beqq}{\begin{eqnarray*}}
\newcommand{\eeq}{\end{eqnarray}}
\newcommand{\eeqq}{\end{eqnarray*}}

\newcommand{\ds}{\displaystyle}

\newcounter{minutes}
\divide\time by 60
\newcounter{hours}
\multiply\time by 60 \addtocounter{minutes}{-\time}

\begin{document}

\bibliographystyle{amsplain}
\title []
{The Schwarz type Lemmas  and the Landau type theorem of mappings
satisfying Poisson's equations}

\def\thefootnote{}
\footnotetext{ \texttt{\tiny File:~\jobname .tex,
          printed: \number\day-\number\month-\number\year,
          \thehours.\ifnum\theminutes<10{0}\fi\theminutes}
} \makeatletter\def\thefootnote{\@arabic\c@footnote}\makeatother

\author{Shaolin Chen}
 \address{Sh. Chen, College of Mathematics and
Statistics, Hengyang Normal University, Hengyang, Hunan 421008,
People's Republic of China.} \email{mathechen@126.com}

\author{David Kalaj}
\address{D. Kalaj, Faculty of Natural Sciences and Mathematics,
University of Montenegro, Cetinjski put b. b. 81000 Podgorica,
Montenegro. } \email{davidk@t-com.me}





\begin{abstract}
For a  given continuous function $g:~\Omega\rightarrow\mathbb{C}$,
we establish some Schwarz type Lemmas for   mappings $f$ in $\Omega$
satisfying the {\rm PDE}: $\Delta f=g$, where $\Omega$ is a subset
of the complex plane $\mathbb{C}$. Then we apply these
 results to obtain a  Landau type theorem, which is  a partial answer to the open problem in
\cite{CP}. 
\end{abstract}


\maketitle \pagestyle{myheadings} \markboth{ Sh. Chen and D.
Kalaj}{The Schwarz type Lemmas  and the Landau type theorem}

\section{Preliminaries and  main results }\label{csw-sec1}
Let $\mathbb{C}  \cong \mathbb{R}^{2}$ be the complex plane. For
$a\in\mathbb{C}$ and  $r>0$, we let $\ID(a,r)=\{z:\, |z-a|<r\}$ so
that $\mathbb{D}_r:=\mathbb{D}(0,r)$ and thus, $\mathbb{D}:=\ID_1$
denotes the open unit disk in the complex plane $\mathbb{C}$. Let
$\mathbb{T}=\partial\mathbb{D}$ be the boundary of $\mathbb{D}$. We
denote by $\mathcal{C}^{m}(\Omega)$ the set of all complex-valued
$m$-times continuously differentiable functions from $\Omega$ into
$\mathbb{C}$, where $\Omega$ is a subset of $\mathbb{C}$  and
$m\in\mathbb{N}_0:=\mathbb{N}\cup\{0\}$. In particular, let
$\mathcal{C}(\Omega):=\mathcal{C}^{0}(\Omega)$, the set of all
continuous functions defined in $\Omega$.

For a real $2\times2$ matrix $A$, we use the matrix norm
$\|A\|=\sup\{|Az|:\,|z|=1\}$ and the matrix function
$\lambda(A)=\inf\{|Az|:\,|z|=1\}$. For $z=x+iy\in\mathbb{C}$, the
formal derivative of the complex-valued functions $f=u+iv$ is given
by
$$D_{f}=\left(\begin{array}{cccc}
\ds u_{x}\;~~ u_{y}\\[2mm]
\ds v_{x}\;~~ v_{y}
\end{array}\right),
$$
so that
$$\|D_{f}\|=|f_{z}|+|f_{\overline{z}}| ~\mbox{ and }~ \lambda(D_{f})=\big| |f_{z}|-|f_{\overline{z}}|\big |,
$$
where $$f_{z}=\frac{\partial f}{\partial z}=\frac{1}{2}\big(
f_x-if_y\big)~\mbox{and}~ f_{\overline{z}}=\frac{\partial
f}{\partial \overline{z}}=\frac{1}{2}\big(f_x+if_y\big).$$  We use
$$J_{f}:=\det D_{f} =|f_{z}|^{2}-|f_{\overline{z}}|^{2}
$$
to denote the {\it Jacobian} of $f$ and $$\Delta
f:=\frac{\partial^{2}f}{\partial
x^{2}}+\frac{\partial^{2}f}{\partial y^{2}}=4f_{z \overline{z}}$$ is
the {\it Laplacian} of $f$ .

For $z, w\in\mathbb{D}$ with $z\neq w$ and $|z|+|w|\neq0$, let
$$G(z,w)=\log\left|\frac{1-z\overline{w}}{z-w}\right|~\mbox{ and}
~P(z,e^{it})=\frac{1-|z|^{2}}{|1-ze^{-it}|^{2}}$$ be the {\it Green
function} and {\it Poisson kernel}, respectively.

Let $\psi:~\mathbb{T}\rightarrow\mathbb{C}$ be a bounded integrable
function  and let $g\in\mathcal{C}(\mathbb{D})$. For
$z\in\mathbb{D}$, the solution to the {\it Poisson's equation}

$$\Delta f(z)=g(z)$$ satisfying the boundary condition
$f|_{\mathbb{T}}=\psi\in L^{1}(\mathbb{T})$ is given by

\be\label{eq-1.0} f(z)=\mathcal{P}_{\psi}(z)-\mathcal{G}_{g}(z),\ee
where
\be\label{eq-2.0}\mathcal{G}_{g}(z)=\frac{1}{2\pi}\int_{\mathbb{D}}G(z,w)g(w)dA(w),~
~\mathcal{P}_{\psi}(z)=\frac{1}{2\pi}\int_{0}^{2\pi}P(z,e^{it})\psi(e^{it})dt,\ee
and $dA(z)$ denotes the Lebesgue measure in $\mathbb{D}$. It is well
known that if $\psi$ and $g$ are continuous in $\mathbb{T}$ and in
$\overline{\mathbb{D}}$, respectively, then
$f=\mathcal{P}_{\psi}-\mathcal{G}_{g}$ has a continuous extension
$\tilde{f}$ to the boundary, and $\tilde{f}=\psi$ in $\mathbb{T}$
(see \cite[pp. 118-120]{Ho}  and \cite{ATM,K2,K1,K4}).


Heinz in his classical paper \cite{He} proved the following result,
which is called the {\it Schwraz Lemma} of complex-valued harmonic
functions: If $f$ is a complex-valued harmonic function from
$\mathbb{D}$ into itself satisfying the condition $f(0)=0$, then,
for $z\in\mathbb{D}$, \be\label{eqh}|f(z)|\leq\frac{4}{\pi}\arctan
|z|.\ee Later, Pavlovi\'c \cite[Theorem 3.6.1]{Pav1} removed the
assumption $f(0)=0$ and improved (\ref{eqh}) into the following
sharp form

\be\label{eq-pav1}\left|f(z)-\frac{1-|z|^{2}}{1+|z|^{2}}f(0)\right|\leq\frac{4}{\pi}\arctan
|z|,\ee where $f$ is a complex-valued harmonic function from
$\mathbb{D}$ into itself.

 The first aim of this paper is to extend (\ref{eq-pav1}) into
mappings satisfying the Poisson's equation, which is as follows.

\begin{thm}\label{thm1.0}
For a given $g\in\mathcal{C}(\overline{\mathbb{D}})$, if
$f\in\mathcal{C}^{2}(\mathbb{D})$ satisfies $\Delta f=g$ and
$f|_{\mathbb{T}}=\psi\in L^{1}(\mathbb{T})$, then, for
$z\in\overline{\mathbb{D}}$,
\be\label{eq-3.0}\left|f(z)-\frac{1}{2\pi}\frac{1-|z|^{2}}{1+|z|^{2}}\int_{0}^{2\pi}\psi(e^{it})dt\right|
\leq\frac{4\|\mathcal{P}_{\psi}\|_{\infty}}{\pi}\arctan|z|+\frac{\|g\|_{\infty}}{4}(1-|z|^{2}),\ee
where
$$\mathcal{P}_{\psi}(z)=\frac{1}{2\pi}\int_{0}^{2\pi}P(z,e^{it})\psi(e^{it})dt,~\|\mathcal{P}_{\psi}\|_{\infty}=
\sup_{z\in\mathbb{D}}|\mathcal{P}_{\psi}(z)|~\mbox{and}~\|g\|_{\infty}=\sup_{z\in\mathbb{D}}|g(z)|.$$
If we take $g(z)=-4M$ and $f(z)=M(1-|z|^{2})$ for
$z\in\overline{\mathbb{D}}$, then the inequality {\rm(\ref{eq-3.0})}
is sharp in $\mathbb{T}\cup\{0\}$, where $M$ is a positive constant.
\end{thm}

The following result is  a classical Schwarz Lemma at the boundary.

\begin{Thm}{\rm (see \cite{G})}\label{Thm-S} Let $f$ be a holomorphic function from
$\mathbb{D}$ into itself. If $f$ is holomorphic at $z=1$ with
$f(0)=0$ and $f(1)=1$, then $f'(1)\geq1$. Moreover, the inequality
is sharp.
\end{Thm}

Theorem \Ref{Thm-S} has attracted much attention and has been
generalized in various forms ( see \cite{BK,Kra,LWT,LT} for
holomorphic functions, and see \cite{K3} for harmonic functions). In
the following, applying Theorem \ref{thm1.0}, we establish a Schwarz
Lemma at the boundary for mappings satisfying the Poisson's
equation, which is  a generalization of Theorem \Ref{Thm-S}.



\begin{thm}\label{thm-h} For a given
$g\in\mathcal{C}(\overline{\mathbb{D}})$, let
$f\in\mathcal{C}^{2}(\mathbb{D})\cap\mathcal{C}(\mathbb{T})$ be a
function of $\mathbb{D}$ into itself satisfying $\Delta f=g,$  where
$\|g\|_{\infty}<\frac{8}{3\pi}.$  If $f(0)=0$ and, for some
$\zeta\in\mathbb{T}$, $\lim_{r\rightarrow1^{-}}|f(r\zeta)|=1$, then

\be\label{eqh5}\liminf_{r\rightarrow1^{-}}\frac{|f(\zeta)-f(r\zeta)|}{1-r}\geq\frac{2}{\pi}-\frac{3\|g\|_{\infty}}{4},
\ee where $r\in[0,1)$.

In particular, if $\|g\|_{\infty}=0$, then the estimate of
{\rm(\ref{eqh5})} is sharp.
\end{thm}

In \cite{Co}, Colonna proved a sharp {\it Schwraz-Pick} type Lemma
of complex-valued harmonic functions, which is as follows: If $f$ is
a complex-valued harmonic function from $\mathbb{D}$ into itself,
then, for $z\in\mathbb{D}$,

\be\label{eq-Co}
\|D_{f}(z)\|\leq\frac{4}{\pi}\frac{1}{1-|z|^{2}}.\ee We extend
(\ref{eq-Co}) into the following form.

\begin{thm}\label{thm2.0}
For a given $g\in\mathcal{C}(\overline{\mathbb{D}})$, if
$f\in\mathcal{C}^{2}(\mathbb{D})$ satisfies $\Delta f=g$ and
$f|_{\mathbb{T}}=\psi\in L^{1}(\mathbb{T})$, then, for
$z\in\mathbb{D}\backslash\{0\},$

\be\label{eq-5.0}\|D_{f}(z)\|\leq\frac{4\|\mathcal{P}_{\psi}\|_{\infty}}{\pi}\frac{1}{1-|z|^{2}}+2\mu(|z|),\ee
where
$$\frac{\|g\|_{\infty}}{4}\leq\mu(|z|)=\frac{\|g\|_{\infty}(1-|z|^{2})}{8|z|^{2}}
\bigg[\frac{1+|z|^{2}}{1-|z|^{2}}-\frac{(1-|z|^{2})}{2|z|}\log\frac{1+|z|}{1-|z|}\bigg]\leq\frac{\|g\|_{\infty}}{3}$$
and  $\mu(|z|)$ is decreasing on $|z|\in(0,1)$. In particular, if
$z=0$, then

\be\label{eq-5.0t}\|D_{f}(0)\|\lim_{|z|\rightarrow0^{+}}\bigg(\frac{4\|\mathcal{P}_{\psi}\|_{\infty}}{\pi}\frac{1}{1-|z|^{2}}+2\mu(|z|)\bigg)=
\frac{4}{\pi}\|\mathcal{P}_{\psi}\|_{\infty}+\frac{2}{3}\|g\|_{\infty}.
\ee Moreover, if $\|g\|_{\infty}=0$, then the extreme functions
$$f(z)=\frac{2M\alpha}{\pi}\arg\left(\frac{1+\phi(z)}{1-\phi(z)}\right)$$
show that the estimate of {\rm(\ref{eq-5.0})} and
{\rm(\ref{eq-5.0t})} are sharp, where $|\alpha|=1$ and $M>0$ are
constants, and $\phi$ is a conformal automorphism of $\mathbb{D}$.
\end{thm}

We remark that if $\|g\|_{\infty}=0$ and
$\|\mathcal{P}_{\psi}\|_{\infty}=1$ in Theorem \ref{thm2.0}, then
(\ref{eq-5.0}) and (\ref{eq-5.0t}) coincide with (\ref{eq-Co}).

Let $\mathcal{A}$  denote the set of all analytic functions $f$
defined in $\mathbb{D}$ satisfying the standard normalization:
$f(0)=f'(0)-1=0$. In the early 20th century, Landau \cite{L} showed
that there is a constant $r>0$, independent of $f\in\mathcal{A}$,
such that $f(\mathbb{D})$ contains a disk of radius $r$. Let $L_{f}$
be the supremum of the set of positive numbers $r$ such that
$f(\mathbb{D})$ contains a disk of radius $r$, where
$f\in\mathcal{A}$. Then we call $\inf_{f\in\mathcal{A}}L_{f}$ the
Landau-Bloch constant. One of the long standing open problems in
geometric function theory is to determine the precise value of the
Landau-Bloch constant. It has attracted much attention, see
\cite{Bo2,LM, M1,M2,W} and references therein.  The Landau theorem
is an important tool in geometric function theory of one complex
variable (cf. \cite{Br,Z}). Unfortunately, for general class of
functions, there is no Landau type theorem (see \cite{HG,W}). In
order to obtain some analogs of the Landau type  theorem for more
general classes of functions, it is necessary to restrict the class
of functions considered (cf.
\cite{AA,Bo1,HG,HG1,CPW-211,CPW0,CV,CP,W}). Let's recall some known
results as follows.

\begin{Thm}{\rm (\cite[Theorem 2]{HG})}\label{Thm-HG}
Let $f$ be a harmonic mapping in $\mathbb{D}$ such that
$f(0)=J_{f}(0)-1=0$ and $|f(z)|<M$ for $z\in\mathbb{D}$, where $M$
is a positive constant. Then $f$ is univalent in
$\mathbb{D}_{\rho_{0}}$ with $\rho_{0}=\pi^{3}/(64mM^{2})$, and
$f(\mathbb{D}_{\rho_{0}})$ contains a univalent disk
$\mathbb{D}_{R_{0}}$ with
$$R_{0}=\frac{\pi}{8M}\rho_{0}=\frac{\pi^{4}}{512mM^{3}},$$ where
$m\approx6.85$ is the minimum of the function
$(3-r^{2})/[r(1-r^{2})]$ for $r\in (0,1)$.
\end{Thm}

\begin{Thm}{\rm (\cite[Theorem 1]{AA})}\label{Thm-AA}
Let $f(z)=|z|^{2}G(z)+K(z)$ be a biharmonic mapping, that is
$\Delta(\Delta f)=0$, in $\mathbb{D}$ such that $f(0)=
K(0)=J_{f}(0)-1=0$, where $G$ and $K$ are harmonic satisfying
$|G(z)|, ~|K(z)|<M$ for $z\in\mathbb{D}$, where $M$ is a positive
constant. Then there is a constant $\rho_{2}\in(0,1)$ so that $f$ is
univalent in $\mathbb{D}_{\rho_{2}}$. In specific $\rho_{2}$
satisfies
$$\frac{\pi}{4M}-2\rho_{2}M-2M\left[\frac{\rho_{2}^{2}}
{(1-\rho_{2})^{2}}+\frac{1}{(1-\rho_{2})^{2}}-1\right]=0$$ and
$f(\mathbb{D}_{\rho_{2}})$ contains a disk $\mathbb{D}_{R_{2}}$,
where
$$R_{2}=\frac{\pi}{4M}\rho_{2}-2M\frac{\rho_{2}^{3}+\rho_{2}^{2}}{1-\rho_{2}}.$$
\end{Thm}

For some $g\in\mathcal{C}(\overline{\mathbb{D}})$, let
$\mathcal{F}_{g}(\overline{\mathbb{D}})$ denote the class of all
complex-valued functions
$f\in\mathcal{C}^{2}(\mathbb{D})\cap\mathcal{C}(\mathbb{T})$
satisfying $\Delta f=g$ and $f(0)=J_{f}(0)-1=0$. We extend Theorems
\Ref{Thm-HG} and   \Ref{Thm-AA} into the following from.

\begin{thm}\label{thm-1}
For a given $g\in\mathcal{C}(\overline{\mathbb{D}})$, let
$f\in\mathcal{F}_{g}(\overline{\mathbb{D}})$ satisfying
$\|g\|_{\infty}\leq M_{1}$ and $\|f\|_{\infty}\leq M_{2}$, where
$M_{1}\geq0$ and $M_{2}>0$ are constants. Then $f$ is univalent in
$\mathbb{D}_{r_{0}}$, where
  $r_{0}$ satisfies the
following equation
$$\frac{1}{\frac{4}{\pi}M_{2}+\frac{2}{3}M_{1}}-\frac{4M_{2}}{\pi}\frac{r_{0}(2-r_{0})}{(1-r_{0})^{2}}
-2M_{1}\big[\log4(1+r_{0})- \log r_{0}\big](2+r_{0})r_{0}=0.$$
Moreover, $f(\mathbb{D}_{r_{0}})$ contains an univalent disk
$\mathbb{D}_{R_{0}}$ with
$$R_{0}\geq\frac{2M_{2}}{\pi}\frac{r_{0}^{2}(2-r_{0})}{(1-r_{0})^{2}}.$$
\end{thm}

\begin{rem}
 Theorem \ref{thm-1} gives an
affirmative answer to the open problem of \cite{CP} for the {\it
$u$-gradient mapping} $f\in\mathcal{C}^{2}(\mathbb{D})$ . If $g$ is
harmonic, then all $f\in\mathcal{F}_{g}(\overline{\mathbb{D}})$ are
biharmonic. Furthermore, if $\|g\|_{\infty}=0$, then all
$f\in\mathcal{F}_{g}(\overline{\mathbb{D}})$ are harmonic. Hence,
Theorem \ref{thm-1}  is also a generalization of series known
results, such as  \cite[Theorem 2]{AA},
 \cite[Theorems, 3, 4, 5 and 6]{HG},
\cite[Theorems 2 and 3]{HG1}, and so on.
\end{rem}

We want to point out that it is failure of Landau type Theorem for
$f\in\mathcal{F}_{g}(\overline{\mathbb{D}})$ without any other
additional condition. It means that the condition
$f(0)=J_{f}(0)-1=0$ is not sufficient to ensure the function of $f$
to the Poisson equation with Landau type theorem, even when
$\|g\|_{\infty}=0$. In particular, Gauthier and Pouryayevali
\cite{GP} proved that it is also failure of Landau's theorem for
quasiconformal mappings $f$ defined in $\mathbb{D}$ satisfying
$f(0)=J_{f}(0)-1=0$.

\begin{exam}
For $g\equiv1$ and $z=x+iy\in\mathbb{D},$ let
$f_{k}(z)=kx+|z|^{2}/4+i\frac{y}{k},$ where $k\in\{1,2,\ldots\}$.
Then, for all $k\in\{1,2,\ldots\}$, $f_{k}$ is univalent. For all
$k\in\{1,2,\ldots\}$, by simple calculations, we see that $J_{f_{k}}
(0)-1=f_{k}(0)=0$, and there is no an absolute constant $\rho_{0}>0$
such that $\mathbb{D}_{\rho_{0}} $ belongs to $f_{k}(\mathbb{D})$.

\end{exam}

\begin{exam}
For $\|g\|_{\infty}=0$ and $z=x+iy\in\mathbb{D},$ let
$f_{k}(z)=kx+i\frac{y}{k},$ where $k\in\{1,2,\ldots\}$. For all
$k\in\{1,2,\ldots\}$, it is not difficult to see that $f_{k}$ is
univalent and  $J_{f_{k}} (0)-1=f_{k}(0)=0$. Moreover, for all
$k\in\{1,2,\ldots\}$, $f_{k}(\mathbb{D})$ contains no disk with
radius bigger than $1/k$. Hence, for all  $k\in\{1,2,\ldots\}$,
there is no an absolute constant $r_{0}>0$ such that
$\mathbb{D}_{r_{0}} $ belongs to $f_{k}(\mathbb{D})$.

\end{exam}


\begin{cor}\label{cor-1}
Under the same hypothesis of Theorem {\rm\ref{thm-1}},  there is a
$r_{0}\in(0,1)$ such that $f$ is bi-Lipschitz in
$\mathbb{D}_{r_{0}}$.
\end{cor}

The proofs of Theorems \ref{thm1.0}, \ref{thm-h}, \ref{thm2.0},
\ref{thm-1} and Corollary \ref{cor-1} will be presented in Section
\ref{csw-sec2}.


\section{Proofs of the  main results }\label{csw-sec2}

\subsection*{Proof of Theorem \ref{thm1.0}} For a given $g\in\mathcal{C}(\mathbb{D})$, by
(\ref{eq-1.0}), we have

\be\label{eq-ck-1}f(z)=\mathcal{P}_{\psi}(z)-\mathcal{G}_{g}(z),
~z\in\mathbb{D},\ee where $\mathcal{P}_{\psi}$ and $\mathcal{G}_{g}$
are defined  in (\ref{eq-2.0}).
  Since $\mathcal{P}_{\psi}$ is
harmonic in $\mathbb{D}$, by (\ref{eq-pav1}), we see that, for
$z\in\mathbb{D}$,

\be\label{eq-y1}
\left|\mathcal{P}_{\psi}(z)-\frac{1-|z|^{2}}{1+|z|^{2}}\mathcal{P}_{\psi}(0)\right|\leq\frac{4\|\mathcal{P}_{\psi}\|_{\infty}}{\pi}\arctan|z|.\ee

On the other hand, for a fixed $z\in\mathbb{D}$, let
$$\zeta=\frac{z-w}{1-\overline{z}w},$$ which is equivalent to $$w=\frac{z-\zeta}{1-\overline{z}\zeta}.$$ Then

\beq\label{eq-y2}\nonumber
\big|\mathcal{G}_{g}(z)\big|&=&\bigg|\frac{1}{2\pi}\int_{\mathbb{D}}\left(\log\frac{1}{|\zeta|}\right)g\left(\frac{z-\zeta}{1-\overline{z}\zeta}\right)
\frac{(1-|z|^{2})^{2}}{|1-\overline{z}\zeta|^{4}}dA(\zeta)\bigg|\\
\nonumber
&\leq&\frac{\|g\|_{\infty}}{2\pi}\bigg|\int_{\mathbb{D}}\left(\log\frac{1}{|\zeta|}\right)
\frac{(1-|z|^{2})^{2}}{|1-\overline{z}\zeta|^{4}}dA(\zeta)\bigg|\\
\nonumber
&=&(1-|z|^{2})^{2}\|g\|_{\infty}\int_{0}^{1}\left[\bigg(\frac{1}{2\pi}\int_{0}^{2\pi}\frac{dt}{|1-\overline{z}re^{it}|^{4}}\bigg)r\log\frac{1}{r}\right]dr\\
\nonumber &=&
(1-|z|^{2})^{2}\|g\|_{\infty}\int_{0}^{1}\left[\bigg(\frac{1}{2\pi}\int_{0}^{2\pi}\frac{dt}{|(1-\overline{z}re^{it})^{2}|^{2}}\bigg)r\log\frac{1}{r}\right]dr\\
\nonumber &=&
(1-|z|^{2})^{2}\|g\|_{\infty}\int_{0}^{1}\left[\bigg(\frac{1}{2\pi}\int_{0}^{2\pi}\Big|\sum_{n=0}^{\infty}(n+1)(r\overline{z})^{n}e^{int}\Big|^{2}dt\bigg)r\log\frac{1}{r}\right]dr
\\
 \nonumber &=&
(1-|z|^{2})^{2}\|g\|_{\infty}\int_{0}^{1}\bigg(r\log\frac{1}{r}\bigg)\sum_{n=0}^{\infty}(n+1)^{2}|z|^{2n}r^{2n}dr\\
\nonumber
&=&(1-|z|^{2})^{2}\|g\|_{\infty}\sum_{n=0}^{\infty}(n+1)^{2}|z|^{2n}\int_{0}^{1}r^{2n+1}\bigg(\log\frac{1}{r}\bigg)dr
\\
\nonumber
&=&\frac{(1-|z|^{2})^{2}\|g\|_{\infty}}{4}\sum_{n=0}^{\infty}|z|^{2n}\\
&=&\frac{\|g\|_{\infty}}{4}(1-|z|^{2}).
 \eeq
Hence, by (\ref{eq-y1}) and (\ref{eq-y2}), we conclude that

\begin{eqnarray*}
\left|f(z)-\frac{1-|z|^{2}}{1+|z|^{2}}\mathcal{P}_{\psi}(0)\right|&\leq&\left|\mathcal{P}_{\psi}(z)-\frac{1-|z|^{2}}{1+|z|^{2}}\mathcal{P}_{\psi}(0)\right|+\big|\mathcal{G}_{g}(z)\big|\\
&\leq&\frac{4\|\mathcal{P}_{\psi}\|_{\infty}}{\pi}\arctan|z|+\frac{\|g\|_{\infty}}{4}(1-|z|^{2}).
\end{eqnarray*}

Now we prove the sharpness part.  For $z\in\overline{\mathbb{D}}$,
let $$g(z)=-4M~\mbox{and}~f(z)=M(1-|z|^{2}),$$ where $M$ is a
positive constant. Then

\begin{eqnarray*}
\left|f(0)-\frac{1}{2\pi}\int_{0}^{2\pi}\psi(e^{it})dt\right|
&=&|\mathcal{G}_{g}(0)|=\left|\frac{1}{2\pi}\int_{\mathbb{D}}\left(\log\frac{1}{|w|}\right)g(w)dA(w)\right|\\
&=&\frac{2M}{\pi}\int_{0}^{2\pi}dt\int_{0}^{1}r\log\frac{1}{r}dr\\
&=&M\\
&=&\frac{\|g\|_{\infty}}{4},
\end{eqnarray*}
which shows (\ref{eq-3.0}) is  sharp at $z=0$. For $z\in\mathbb{T}$,
the optimality of (\ref{eq-3.0})  is obvious.
 The proof of this theorem is
complete.  \qed

\subsection*{Proof of Theorem \ref{thm-h}} For a given $g\in\mathcal{C}(\overline{\mathbb{D}})$, by
(\ref{eq-1.0}) with $f$ in place of $\psi$, we have

$$f(z)=\mathcal{P}_{f}(z)-\mathcal{G}_{g}(z), ~z\in\mathbb{D},$$ where $\mathcal{P}_{f}$
 and $\mathcal{G}_{g}$ are defined
in (\ref{eq-2.0}). Since $f(0)=0$, we see that

\beq\label{eq-h2}|\mathcal{P}_{f}(0)|&=&|\mathcal{G}_{g}(0)|=\left|\frac{1}{2\pi}\int_{\mathbb{D}}\log\frac{1}{|w|}g(w)dA(w)\right|\\
\nonumber
&\leq&\frac{\|g\|_{\infty}}{2\pi}\int_{0}^{2\pi}dt\int_{0}^{1}r\log\frac{1}{r}dr\\
\nonumber &=&\frac{\|g\|_{\infty}}{4}. \eeq

By (\ref{eq-h2}) and Theorem \ref{thm1.0}, we have

\begin{eqnarray*}
|f(\zeta)-f(r\zeta)|&=&\left|f(\zeta)+\mathcal{P}_{f}(0)\frac{1-|z|^{2}}{1+|z|^{2}}-\mathcal{G}_{g}(0)\frac{1-|z|^{2}}{1+|z|^{2}}-f(r\zeta)\right|\\
&\geq&1-\left|f(r\zeta)-\mathcal{P}_{f}(0)\frac{1-|z|^{2}}{1+|z|^{2}}\right|-|\mathcal{G}_{g}(0)|\frac{1-|z|^{2}}{1+|z|^{2}}\\
&\geq&1-\frac{4}{\pi}\arctan|z|-\frac{\|g\|_{\infty}}{4}(1-|z|^{2})-|\mathcal{G}_{g}(0)|\frac{1-|z|^{2}}{1+|z|^{2}}\\
&\geq&1-\frac{4}{\pi}\arctan|z|-\frac{\|g\|_{\infty}}{4}(1-|z|^{2})-\frac{\|g\|_{\infty}}{4}\frac{(1-|z|^{2})}{1+|z|^{2}},
\end{eqnarray*}

which, together with L'Hopital's rule, gives that

\begin{eqnarray*}
\liminf_{r\rightarrow1^{-}}\frac{|f(e^{i\theta})-f(re^{i\theta})|}{1-r}
&\geq&\lim_{r\rightarrow1^{-}}\frac{1-\frac{4}{\pi}\arctan
r-\frac{\|g\|_{\infty}}{4}(1-r^{2})-\frac{\|g\|_{\infty}}{4}\frac{(1-r^{2})}{1+r^{2}}}{1-r}\\
&=&\lim_{r\rightarrow1^{-}}\left[\frac{4}{\pi}\frac{1}{1+r^{2}}-\|g\|_{\infty}\frac{r}{2}-\|g\|_{\infty}\frac{r}{(1+r^{2})^{2}}\right]\\
&=&\frac{2}{\pi}-\frac{\|g\|_{\infty}}{2}-\frac{\|g\|_{\infty}}{4}\\
&=&\frac{2}{\pi}-\frac{3\|g\|_{\infty}}{4},
\end{eqnarray*}
where $z=r\zeta$ and $\zeta\in\mathbb{T}$.

The sharpness part  easily follows from \cite[Theorem 2.5]{K3}. The
proof of this theorem is complete. The proof of this theorem is
complete.
 \qed

\begin{Thm}{\rm (\cite{T}~$\mbox{or}$~\cite[Proposition 2.4]{K2})}\label{Thm-B}
Let $X$ be an open subset of $\mathbb{R}$, and $\Omega$ be a measure
space. Suppose that a function
$F:~X\times\Omega\rightarrow\mathbb{R}$ satisfies the following
conditions:

\begin{enumerate}
\item[(1)] $F(x,w)$ is a measurable function of $x$ and $w$ jointly,
and is integrable over $\omega$, for almost all $x\in X$ held fixed.
\item[(2)] For almost all $w\in\Omega$, $F(x,w)$ is an absolutely
continuous function of $x$. $($This guarantees that $\partial
F(x,w)/\partial x$ exists almost everywhere.$)$
\item[(3)] $\partial F/\partial x$ is locally integrable; that is,
for all compact intervals $[a,b]$ contained in $X$:

$$\int_{a}^{b}\int_{\Omega}\left|\frac{\partial}{\partial x}F(x,w)\right|dwdx<\infty.$$
\end{enumerate}
Then $\int_{\Omega}F(x,w)dw$ is an absolutely continuous function of
$x$, and for almost every $x\in X$, its derivative exists and is
given by
$$\frac{d}{dx}\int_{\Omega}F(x,w)dw=\int_{\Omega}\frac{\partial}{\partial x}F(x,w)dw.$$
\end{Thm}


\subsection*{Proof of Theorem \ref{thm2.0}}
For a given $g\in\mathcal{C}(\overline{\mathbb{D}})$, by
(\ref{eq-ck-1}), we have

$$f(z)=\mathcal{P}_{\psi}(z)-\mathcal{G}_{g}(z), ~z\in\mathbb{D},$$ where
$\mathcal{P}_{\psi}$ and $\mathcal{G}_{g}$ are the same as in
(\ref{eq-ck-1}). Applying \cite[Lemma 2.3]{K2} and Theorem
\Ref{Thm-B}, we have

\beq\label{eqcx1} \frac{\partial}{\partial
z}\mathcal{G}_{g}(z)&=&\frac{1}{2\pi}\int_{\mathbb{D}}\frac{\partial}{\partial
z}G(z,w)g(w)dA(w)\\ \nonumber
&=&\frac{1}{4\pi}\int_{\mathbb{D}}\frac{(1-|w|^{2})}{(z-w)(z\overline{w}-1)}g(w)dA(w)\in\mathcal{C}(\mathbb{D})
\eeq and
\begin{eqnarray*}
\frac{\partial}{\partial \overline{z}}\mathcal{G}_{g}(z)&=&
\frac{1}{2\pi}\int_{\mathbb{D}}\frac{\partial}{\partial
\overline{z}}G(z,w)g(w)dA(w)\\
&=&\frac{1}{4\pi}\int_{\mathbb{D}}\frac{(1-|w|^{2})}{(\overline{z}-\overline{w})(w\overline{z}-1)}g(w)dA(w)\in\mathcal{C}(\mathbb{D}).
\end{eqnarray*}

For a fixed $z\in\mathbb{D}\backslash\{0\}$, let
$$\zeta=\frac{z-w}{1-\overline{z}w}$$ which implies that

\be\label{eq-h4c}w=\frac{z-\zeta}{1-\overline{z}\zeta},~1-\overline{z}w=\frac{1-|z|^{2}}{1-\overline{z}\zeta}~
\mbox{and}~1-|w|^{2}=\frac{(1-|\zeta|^{2})(1-|z|^{2})}{|1-\overline{z}\zeta|^{2}}.\ee

Then, by (\ref{eqcx1}), (\ref{eq-h4c}) and the change of variables,
we have

\beq\label{eqcx2} \nonumber\left|\frac{\partial}{\partial
z}\mathcal{G}_{g}(z)\right|&\leq&\frac{1}{4\pi}\int_{\mathbb{D}}\frac{(1-|w|^{2})}{|z-w||z\overline{w}-1|}|g(w)|dA(w)\\
\nonumber
&\leq&\frac{\|g\|_{\infty}}{4\pi}\int_{\mathbb{D}}\frac{(1-|w|^{2})}{|z-w||z\overline{w}-1|}dA(w)\\
\nonumber
&=&\frac{\|g\|_{\infty}}{4\pi}\int_{\mathbb{D}}\frac{(1-|w|^{2})}{|\zeta||1-\overline{z}w|^{2}}\frac{(1-|z|^{2})^{2}}{|1-\overline{z}\zeta|^{4}}dA(\zeta)
\\ \nonumber
&=&\frac{\|g\|_{\infty}}{4\pi}\int_{\mathbb{D}}\frac{(1-|z|^{2})(1-|\zeta|^{2})}{|\zeta||1-\overline{z}\zeta|^{4}}dA(\zeta)\\
\nonumber
&=&\frac{\|g\|_{\infty}(1-|z|^{2})}{2}\int_{0}^{1}\left[(1-r^{2})\bigg(\frac{1}{2\pi}\int_{0}^{2\pi}\frac{dt}{|1-\overline{z}re^{it}|}\bigg)\right]dr\\
\nonumber &=&\frac{\|g\|_{\infty}(1-|z|^{2})}{2}\int_{0}^{1}
\left[(1-r^{2})\bigg(\frac{1}{2\pi}\int_{0}^{2\pi}\Big|\sum_{n=0}^{\infty}(n+1)(r\overline{z})^{n}e^{int}\Big|^{2}dt\bigg)\right]dr\\
\nonumber &=&\frac{\|g\|_{\infty}(1-|z|^{2})}{2}\int_{0}^{1}
(1-r^{2})\left[\sum_{n=0}^{\infty}(n+1)^{2}|z|^{2n}r^{2n}\right]dr\\
\nonumber &=&\frac{\|g\|_{\infty}(1-|z|^{2})}{2}\int_{0}^{1}
\frac{(1-r^{2})(1+|z|^{2}r^{2})}{(1-|z|^{2}r^{2})^{3}}dr\\
&=&\frac{\|g\|_{\infty}(1-|z|^{2})}{2}\bigg[-\frac{1}{|z|^{2}}I_{1}+
\left(\frac{3}{|z|^{2}}-1\right)I_{2}
+2\left(1-\frac{1}{|z|^{2}}\right)I_{3}\bigg], \eeq where


\beq\label{eqcx3}
I_{1}=\int_{0}^{1}\frac{dr}{1-r^{2}|z|^{2}}=\frac{1}{|z|}\log\frac{1+|z|r}{\sqrt{1-|z|^{2}r^{2}}}\bigg|_{0}^{1}=\frac{1}{|z|}\log\frac{1+|z|}{\sqrt{1-|z|^{2}}},\eeq

\beq\label{eqcx4} \nonumber
I_{2}=\int_{0}^{1}\frac{dr}{(1-r^{2}|z|^{2})^{2}}&=&\frac{1}{2|z|}\left(\log\frac{1+|z|r}{\sqrt{1-|z|^{2}r^{2}}}+
\frac{|z|r}{1-|z|^{2}r^{2}}\right)\bigg|_{0}^{1}\\
&=&\frac{1}{2|z|}\log\frac{1+|z|}{\sqrt{1-|z|^{2}}}+\frac{1}{2(1-|z|^{2})}\eeq
and

\beq\label{eqcx5}\nonumber
I_{3}=\int_{0}^{1}\frac{dr}{(1-r^{2}|z|^{2})^{3}}&=&\frac{1}{4|z|}\left(\frac{|z|r}{(1-r^{2}|z|^{2})^{2}}
+\frac{3}{2}\frac{|z|r}{1-r^{2}|z|^{2}}+\frac{3}{2}\log\frac{1+|z|r}{\sqrt{1-|z|^{2}r^{2}}}\right)\bigg|_{0}^{1}\\
&=&\frac{1}{4(1-|z|^{2})^{2}}+\frac{3}{8(1-|z|^{2})}+\frac{3}{8|z|}\log\frac{1+|z|}{\sqrt{1-|z|^{2}}}.
\eeq

By (\ref{eqcx3}), (\ref{eqcx4}) and (\ref{eqcx5}), we get

\beq\nonumber-\frac{1}{|z|^{2}}I_{1}+
\left(\frac{3}{|z|^{2}}-1\right)I_{2}
+2\left(1-\frac{1}{|z|^{2}}\right)I_{3}&=&\frac{1}{4|z|^{2}}\bigg[\frac{1+|z|^{2}}{1-|z|^{2}}\\
\nonumber
&&-\frac{(1-|z|^{2})}{2|z|}\log\frac{1+|z|}{1-|z|}\bigg],\eeq which,
together with (\ref{eqcx2}), yields that

\be\label{eqcx6}\left|\frac{\partial}{\partial
z}\mathcal{G}_{g}(z)\right|\leq\mu(|z|),\ee where
$$\mu(|z|)=\frac{\|g\|_{\infty}(1-|z|^{2})}{8|z|^{2}}
\bigg[\frac{1+|z|^{2}}{1-|z|^{2}}-\frac{(1-|z|^{2})}{2|z|}\log\frac{1+|z|}{1-|z|}\bigg].$$
 By a similar proof process of (\ref{eqcx6}), we have

\be\label{eqcx7}\left|\frac{\partial}{\partial
\overline{z}}\mathcal{G}_{g}(z)\right|\leq\mu(|z|).\ee

By  direct calculation (or by \cite[Lemma 2.3]{K2}), we obtain

\be\label{eqcx8}\lim_{|z|\rightarrow0^{+}}\frac{\|g\|_{\infty}(1-|z|^{2})}{8|z|^{2}}
\bigg[\frac{1+|z|^{2}}{1-|z|^{2}}-\frac{(1-|z|^{2})}{2|z|}\log\frac{1+|z|}{1-|z|}\bigg]=\frac{\|g\|_{\infty}}{3},\ee

$$\lim_{|z|\rightarrow1^{-}}\frac{\|g\|_{\infty}(1-|z|^{2})}{8|z|^{2}}
\bigg[\frac{1+|z|^{2}}{1-|z|^{2}}-\frac{(1-|z|^{2})}{2|z|}\log\frac{1+|z|}{1-|z|}\bigg]=\frac{\|g\|_{\infty}}{4}$$
and  $\mu(|z|)$ is decreasing on $|z|\in(0,1)$.

On the other hand, since $\mathcal{P}_{\psi}$ is harmonic in
$\mathbb{D}$, by \cite[Theorem 3]{Co} (see also \cite{CV,CPRW}), we
see that, for $z\in\mathbb{D}$,
\be\label{eq-7.0}\|D_{\mathcal{P}_{\psi}}(z)\|\leq\frac{4\|\mathcal{P}_{\psi}\|_{\infty}}{\pi}\frac{1}{1-|z|^{2}}.\ee
Hence (\ref{eq-5.0}) follows from (\ref{eqcx6}), (\ref{eqcx7}) and
(\ref{eq-7.0}). Furthermore, applying (\ref{eq-5.0}) and
(\ref{eqcx8}), we get (\ref{eq-5.0t}).   The proof of this theorem
is complete. \qed

New wi formulate the following well-known result
\begin{lem}\label{lem-0} The improper integral
$$\int_{0}^{\frac{\pi}{2}}\log\sin x
dx=\int_{0}^{\frac{\pi}{2}}\log\cos x dx=-\frac{\pi}{2}\log2.$$
\end{lem}

\begin{lem}\label{lem-1}
For $z\in\mathbb{D}\backslash\{0\}$, the improper integral
\begin{eqnarray*}
\int_{\mathbb{D}}\frac{dA(w)}{|w||z-w|}&=&\int_{0}^{2\pi}\log\big(1-r\cos
t+\sqrt{1+r^{2}-2r\cos t}\big)dt\\
&& -2\pi \log r+2\pi\log2\\
&\leq&2\pi\log4(1+r)-2\pi \log r,
\end{eqnarray*} where
$r=|z|$.
\end{lem}

\bpf Let $z=re^{i\alpha}$ and $w=\rho e^{i\theta} $. Then
\beq\label{eq-2}\nonumber
\int_{\mathbb{D}}\frac{dA(w)}{|w||z-w|}&=&\int_{0}^{1}d\rho\int_{0}^{2\pi}\frac{d\theta}{\sqrt{r^{2}+\rho^{2}-2\rho
r\cos (\theta-\alpha)}}\\ \nonumber
&=&\int_{0}^{1}d\rho\int_{0}^{2\pi}\frac{dt}{\sqrt{r^{2}+\rho^{2}-2\rho
r\cos t}}\\
\nonumber
&=&\int_{0}^{2\pi}dt\int_{0}^{1}\frac{d\rho}{\sqrt{r^{2}+\rho^{2}-2\rho
r\cos t}}\\
\nonumber &=&\int_{0}^{2\pi}\bigg\{\frac{1}{2r \cos
t}\bigg[\int_{0}^{1}\frac{2\rho d\rho}{\sqrt{r^{2}+\rho^{2}-2\rho
r\cos
t}}\\
\nonumber&&-\int_{0}^{1}\frac{d(r^{2}+\rho^{2}-2\rho r\cos
t)}{\sqrt{r^{2}+\rho^{2}-2\rho r\cos t}}\bigg]\bigg\}dt\\
\nonumber&=&\int_{0}^{2\pi}\bigg[\frac{1}{r\cos
t}\int_{0}^{1}\frac{\rho d\rho}{\sqrt{r^{2}+\rho^{2}-2\rho r\cos t}}
\\
\nonumber&&-\frac{1}{r\cos t}\big(\sqrt{1+r^{2}-2r\cos
t}-r\big)\bigg]dt
\\
&=&\int_{0}^{2\pi}\bigg[\frac{1}{r\cos t}\int_{0}^{1}\frac{\rho
d\rho}{\sqrt{r^{2}+\rho^{2}-2\rho r\cos t}}
\\ \nonumber
&&-\frac{\sqrt{1+r^{2}-2r\cos t}}{r\cos t}+\frac{1}{\cos t}\bigg]dt.
 \eeq
By calculations, we get
\beq\label{eq-3} \nonumber\int_{0}^{1}\frac{\rho
d\rho}{\sqrt{r^{2}+\rho^{2}-2\rho r\cos t}}&=&H(\rho)|_{0}^{1}\\
 \nonumber &=&\sqrt{1+r^{2}-2r\cos t}\\
\nonumber &&+r\cos t \log\left(1-r\cos t+\sqrt{1+r^{2}-2r\cos
t}\right)\\  &&- r-r\cos t\log r(1-\cos t),
\eeq
 where
$$H(\rho)=\sqrt{\rho^{2}+r^{2}-2r\rho\cos t}+r \cos
t\log\left(\rho-r\cos t+\sqrt{r^{2}+\rho^{2}-2\rho r\cos
t}\right).$$

By (\ref{eq-2}), (\ref{eq-3}) and Lemma \ref{lem-0}, we see that

\beq\label{eq-4}
\int_{\mathbb{D}}\frac{dA(w)}{|w||z-w|}&=&\int_{0}^{2\pi}\log\big(1-r\cos
t+\sqrt{1+r^{2}-2r\cos t}\big)dt\\ \nonumber&&-\int_{0}^{2\pi}\log
r(1-\cos t)dt\\ \nonumber&=&\int_{0}^{2\pi}\log\big(1-r\cos
t+\sqrt{1+r^{2}-2r\cos t}\big)dt
\\ \nonumber&&-2\pi\log r-\int_{0}^{2\pi}\log\left(2\sin^{2}\frac{t}{2}\right)dt
\\ \nonumber&=&
\int_{0}^{2\pi}\log\big(1-r\cos t+\sqrt{1+r^{2}-2r\cos t}\big)dt\\
\nonumber&&-2\pi \log2r-8\int_{0}^{\frac{\pi}{2}}\log(\sin t) dt
\\ \nonumber&=&\int_{0}^{2\pi}\log\big(1-r\cos t+\sqrt{1+r^{2}-2r\cos t}\big)dt
\\
\nonumber&&-2\pi \log r+2\pi\log2\\
\nonumber&\leq&2\pi\log4(1+r)-2\pi \log r.\eeq The proof of this
theorem is complete. \epf

\begin{Lem}{\rm (\cite[Lemma 1]{CPW0})}\label{LemA}
Let $f$ be a harmonic mapping of $\mathbb{D}$ into $\mathbb{C}$ such
that $|f(z)|\leq M$ and
$f(z)=\sum_{n=0}^{\infty}a_{n}z^{n}+\sum_{n=1}^{\infty}\overline{b}_{n}\overline{z}^{n}$.
Then $|a_{0}|\leq M$ and for all $n\geq 1,$
$$|a_{n}|+|b_{n}|\leq \frac{4M}{\pi}.
$$
\end{Lem}

\begin{lem}\label{lem3}
For $x\in(0,1)$, let
$$\phi(x)=\frac{1}{\frac{4}{\pi}M_{2}+\frac{2}{3}M_{1}}-\frac{4M_{2}}{\pi}\frac{x(2-x)}{(1-x)^{2}}
-2M_{1}\big[\log4(1+x)- \log x\big](2+x)x,$$ where $M_{2}>0$ and
$M_{1}\geq0$ are constant. Then $\phi$ is strictly decreasing and
there is an unique $x_{0}\in(0,1)$ such that $\phi(x_{0})=0.$
\end{lem}

\bpf For $x\in(0,1)$, let
$$f_{1}(x)=\frac{4M_{2}}{\pi}\frac{x(2-x)}{(1-x)^{2}}$$ and $$f_{2}(x)=2M_{1}\big[\log2(1+x)- \log x+\log2\big](2+x)x.$$
Since, for $x\in(0,1)$,
$$f_{1}'(x)=\frac{8M_{2}}{\pi}\frac{1}{(1-x)^{3}}>0$$ and

\begin{eqnarray*}
f_{2}'(x)&=&2M_{1}\left[2(x+1)\log\frac{4(1+x)}{x}-\frac{2+x}{1+x}\right]\\
&=&2M_{1}\left\{2(x+1)\bigg[\log4+\log\Big(1+\frac{1}{x}\Big)\bigg]-\frac{2+x}{1+x}\right\}\\
&\geq&2M_{1}\left\{2(x+1)\bigg[1+\frac{1}{1+x}\bigg]-\frac{2+x}{1+x}\right\}\\
&=&2M_{1}\frac{(2+x)(2x+1)}{1+x}\geq0,
\end{eqnarray*}
we see that $f_{1}+f_{2}$ is continuous and strictly increasing in
$(0,1)$. Then $\phi$ is continuous and strictly decreasing in
$(0,1)$, which, together with
$$\lim_{x\rightarrow0^{+}}\phi(x)=\frac{1}{\frac{4}{\pi}M_{2}+\frac{2}{3}M_{1}}~\mbox{and}~\lim_{x\rightarrow1^{-}}\phi(x)=-\infty,$$
implies that there is an unique $x_{0}\in(0,1)$ such that
$\phi(x_{0})=0.$  \epf

\begin{lem}\label{lem4}
For $x\in(0,1]$, let
$$\tau_{1}(x)=\frac{2-r_{0}x}{(1-r_{0}x)^{2}}~\mbox{and}~\tau_{2}(x)=x\big[\log4(1+r_{0}x)-\log
(r_{0}x)\big],$$ where $r_{0}\in(0,1)$ is a constant. Then
$\tau_{1}$ and $\tau_{2}$ are increasing functions in $(0,1]$.
\end{lem}

\subsection*{Proof of Theorem \ref{thm-1}} As before, by
(\ref{eq-ck-1}) with $f$ in place of $\psi$, we have

$$f(z)=\mathcal{P}_{f}(z)-\mathcal{G}_{g}(z), ~z\in\mathbb{D},$$ where
$\mathcal{P}_{f}$ and $\mathcal{G}_{g}$ are defined  in
(\ref{eq-ck-1}).  By \cite[Lemma 2.3]{K2}, Theorem \Ref{Thm-B} and
Lemma \ref{lem-1}, we have

\beq\label{eq-2c} \left|\frac{\partial\mathcal{G}_{g}(z)}{\partial
z}-\frac{\partial\mathcal{G}_{g}(0)}{\partial
z}\right|&=&\bigg|\frac{1}{4\pi}\int_{\mathbb{D}}\frac{(1-|w|^{2})}{(z-w)(z\overline{w}-1)}g(w)dA(w)\\
\nonumber
&&-\frac{1}{4\pi}\int_{\mathbb{D}}\frac{(1-|w|^{2})}{w}g(w)dA(w)\bigg|\\
\nonumber
&=&\bigg|\frac{1}{4\pi}\int_{\mathbb{D}}\frac{z(1-|w|^{2})(1+|w|^{2}-z\overline{w})}{w(z-w)(z\overline{w}-1)}g(w)dA(w)\bigg|\\
\nonumber
&\leq&\frac{M_{1}|z|}{4\pi}\int_{\mathbb{D}}\frac{(1-|w|^{2})\big|1+|w|^{2}-z\overline{w}\big|}{|w||z-w||1-z\overline{w}|}dA(w)\\
\nonumber
&\leq&\frac{|z|(2+|z|)M_{1}}{4\pi}\int_{\mathbb{D}}\frac{(1+|w|)}{|w||z-w|}dA(w)\\
\nonumber
&\leq&\frac{|z|(2+|z|)M_{1}}{2\pi}\int_{\mathbb{D}}\frac{1}{|w||z-w|}dA(w)\\
\nonumber
 &\leq&M_{1}\big[\log4(1+|z|)- \log
|z|\big]|z|(2+|z|). \eeq

By a similar proof process of (\ref{eq-2c}), we get

\beq\label{eq-3c} \left|\frac{\partial\mathcal{G}_{g}(z)}{\partial
\overline{z}}-\frac{\partial\mathcal{G}_{g}(0)}{\partial
\overline{z}}\right|&=&\bigg|\frac{1}{4\pi}\int_{\mathbb{D}}\frac{(1-|w|^{2})}{(\overline{z}-\overline{w})(w\overline{z}-1)}g(w)dA(w)\\
\nonumber
&&-\frac{1}{4\pi}\int_{\mathbb{D}}\frac{(1-|w|^{2})}{\overline{w}}g(w)dA(w)\bigg|\\
\nonumber
&\leq&M_{1}\big[\log4(1+|z|)- \log |z|\big]|z|(2+|z|).
\eeq

On the other hand,  $\mathcal{P}_{f}$ can be written by
$$\mathcal{P}_{f}(z)=\sum_{n=0}^{\infty}a_{n}z^{n}+\sum_{n=1}^{\infty}\overline{b}_{n}\overline{z}^{n}$$ because $\mathcal{P}_{f}$ is harmonic in
$\mathbb{D}$.

Since $|\mathcal{P}_{f}(z)|\leq M_{2}$  for $z\in\mathbb{D}$, by
Lemma \Ref{LemA}, we have

\be\label{eq-5} |a_{n}|+|b_{n}|\leq \frac{4M_{2}}{\pi}\ee for
$n\geq1$.

By (\ref{eq-5}), we see that

\beq\label{eq-6} \nonumber \left|\frac{\partial
\mathcal{P}_{f}(z)}{\partial z}- \frac{\partial
\mathcal{P}_{f}(0)}{\partial z}\right|+\left|\frac{\partial
\mathcal{P}_{f}(z)}{\partial \overline{z}}- \frac{\partial
\mathcal{P}_{f}(0)}{\partial
\overline{z}}\right|&=&\left|\sum_{n=2}^{\infty}na_{n}z^{n-1}\right|+\left|\sum_{n=2}^{\infty}nb_{n}\overline{z}^{n-1}\right|\\
\nonumber
&\leq&\sum_{n=2}^{\infty}n\big(|a_{n}|+|b_{n}|\big)|z|^{n-1}
\\
\nonumber &\leq&\frac{4M_{2}}{\pi}\sum_{n=2}^{\infty}n|z|^{n-1}\\
&=&\frac{4M_{2}}{\pi}\frac{|z|(2-|z|)}{(1-|z|)^{2}}. \eeq

Applying Theorem \ref{thm2.0}, we obtain

$$
1=J_{f}(0)=\|D_{f}(0)\|\lambda(D_{f}(0))\leq
\lambda(D_{f}(0))\left(\frac{4}{\pi}M_{2}+\frac{2}{3}M_{1}\right),$$
which gives that

\be\label{eq-7}
\lambda(D_{f}(0))\geq\frac{1}{\frac{4}{\pi}M_{2}+\frac{2}{3}M_{1}}.\ee

In order to prove the univalence of $f$ in $\mathbb{D}_{r_{0}}$, we
choose two distinct points  $z_{1}, z_{2}\in\mathbb{D}_{r_{0}}$ and
let $[z_{1},z_{2}]$ denote the segment from $z_{1}$ to $z_{2}$ with
the endpoints $z_{1}$ and $z_{2}$, where $r_{0}$ satisfies the
following equation
$$\frac{1}{\frac{4}{\pi}M_{2}+\frac{2}{3}M_{1}}-\frac{4M_{2}}{\pi}\frac{r_{0}(2-r_{0})}{(1-r_{0})^{2}}-2M_{1}\big[\log4(1+r_{0})- \log r_{0}\big](2+r_{0})r_{0}=0.$$ By
(\ref{eq-2c}), (\ref{eq-3c}), (\ref{eq-6}), (\ref{eq-7}), Lemmas
\ref{lem3} and \ref{lem4}, we have


\beq\label{eq-bi-1}
\nonumber|f(z_{2})-f(z_{1})|&=&\left|\int_{[z_{1},z_{2}]}f_{z}(z)dz+f_{\overline{z}}(z)d\overline{z}\right|\\
\nonumber
&=&\left|\int_{[z_{1},z_{2}]}f_{z}(0)dz+f_{\overline{z}}(0)d\overline{z}\right|\\
\nonumber
&&-\left|\int_{[z_{1},z_{2}]}\big(f_{z}(z)-f_{z}(0)\big)dz+\big(f_{\overline{z}}(z)-f_{\overline{z}}(0)\big)d\overline{z}\right|\\
\nonumber &\geq&\lambda(D_{f}(0))|z_{2}-z_{1}|\\ \nonumber
&&-\int_{[z_{1},z_{2}]}\big(|f_{z}(z)-f_{z}(0)|+|f_{\overline{z}}(z)-f_{\overline{z}}(0)|\big)|dz|\\
\nonumber &\geq&\lambda(D_{f}(0))|z_{2}-z_{1}|\\ \nonumber&&
-\int_{[z_{1},z_{2}]}\left(\bigg|\frac{\partial\mathcal{G}_{g}(z)}{\partial
z}-\frac{\partial\mathcal{G}_{g}(0)}{\partial
z}\bigg|+\bigg|\frac{\partial\mathcal{G}_{g}(z)}{\partial
\overline{z}}-\frac{\partial\mathcal{G}_{g}(0)}{\partial
\overline{z}}\bigg|\right)|dz|\\ \nonumber&&
-\int_{[z_{1},z_{2}]}\left(\bigg|\frac{\partial
\mathcal{P}_{f}(z)}{\partial z}- \frac{\partial
\mathcal{P}_{f}(0)}{\partial z}\bigg| +\bigg|\frac{\partial
\mathcal{P}_{f}(z)}{\partial \overline{z}}- \frac{\partial
\mathcal{P}_{f}(0)}{\partial
\overline{z}}\bigg|\right)|dz|\\
&>&|z_{2}-z_{1}|\big\{\lambda(D_{f}(0))-\frac{4M_{2}}{\pi}\frac{r_{0}(2-r_{0})}{(1-r_{0})^{2}}\\
\nonumber &&-2M_{1}\big[\log4(1+r_{0})- \log
r_{0}\big](2+r_{0})r_{0}\big\}\\ \nonumber
&\geq&|z_{2}-z_{1}|\bigg\{\frac{1}{\frac{4}{\pi}M_{2}+\frac{2}{3}M_{1}}-\frac{4M_{2}}{\pi}\frac{r_{0}(2-r_{0})}{(1-r_{0})^{2}}\\
\nonumber &&-2M_{1}\big[\log4(1+r_{0})- \log
r_{0}\big](2+r_{0})r_{0}\bigg\}\\ \nonumber&=&0,
\eeq
which yields that $f(z_{2})\neq f(z_{1})$. The univalence of $f$
follows from the arbitrariness of $z_{1}$ and $z_{2}$.

Now, for all $\zeta=r_{0}e^{i\theta}\in\partial\mathbb{D}_{r_{0}}$,
by (\ref{eq-2c}), (\ref{eq-3c}), (\ref{eq-6}), (\ref{eq-7}), Lemmas
\ref{lem3} and \ref{lem4}, we obtain

\begin{eqnarray*}
|f(\zeta)-f(0)|&=&\left|\int_{[0,\zeta]}f_{z}(z)dz+f_{\overline{z}}(z)d\overline{z}\right|\\
&=&\left|\int_{[0,\zeta]}f_{z}(0)dz+f_{\overline{z}}(0)d\overline{z}\right|\\
&&-\left|\int_{[0,\zeta]}\big(f_{z}(z)-f_{z}(0)\big)dz+\big(f_{\overline{z}}(z)-f_{\overline{z}}(0)\big)d\overline{z}\right|\\
&\geq&\lambda(D_{f}(0))r_{0}\\
&&-\int_{[0,\zeta]}\big(|f_{z}(z)-f_{z}(0)|+|f_{\overline{z}}(z)-f_{\overline{z}}(0)|\big)|dz|\\
&\geq&\lambda(D_{f}(0))r_{0}\\&&
-\int_{[0,\zeta]}\left(\bigg|\frac{\partial\mathcal{G}_{g}(z)}{\partial
z}-\frac{\partial\mathcal{G}_{g}(0)}{\partial
z}\bigg|+\bigg|\frac{\partial\mathcal{G}_{g}(z)}{\partial
\overline{z}}-\frac{\partial\mathcal{G}_{g}(0)}{\partial
\overline{z}}\bigg|\right)|dz|\\&&
-\int_{[0,\zeta]}\left(\bigg|\frac{\partial
\mathcal{P}_{f}(z)}{\partial z}- \frac{\partial
\mathcal{P}_{f}(0)}{\partial z}\bigg| +\bigg|\frac{\partial
\mathcal{P}_{f}(z)}{\partial \overline{z}}- \frac{\partial
\mathcal{P}_{f}(0)}{\partial
\overline{z}}\bigg|\right)|dz|\\
&\geq&\frac{r_{0}}{\frac{4}{\pi}M_{2}+\frac{2}{3}M_{1}}-\frac{4M_{2}}{\pi}\int_{[0,\zeta]}\frac{|z|(2-|z|)}{(1-|z|)^{2}}|dz|\\
&&-2M_{1}\int_{[0,\zeta]}\big[\log4(1+|z|)- \log
|z|\big]|z|(2+|z|)|dz|\\
&=&\frac{r_{0}}{\frac{4}{\pi}M_{2}+\frac{2}{3}M_{1}}-\frac{4M_{2}r_{0}^{2}}{\pi}\int_{0}^{1}\frac{t(2-r_{0}t)}{(1-r_{0}t)^{2}}dt\\
&&-2M_{1}r_{0}^{2}\int_{0}^{1}\big[\log4(1+r_{0}t)- \log
(r_{0}t)\big]t(2+r_{0}t)dt\\
&\geq&\frac{r_{0}}{\frac{4}{\pi}M_{2}+\frac{2}{3}M_{1}}-\frac{4M_{2}r_{0}^{2}}{\pi}\frac{(2-r_{0})}{(1-r_{0})^{2}}\int_{0}^{1}t dt\\
&&-2M_{1}r_{0}^{2}(2+r_{0})\int_{0}^{1}\big[\log4(1+r_{0}t)-
\log (r_{0}t)\big]t dt\\
&\geq&r_{0}\bigg\{\frac{1}{\frac{4}{\pi}M_{2}+\frac{2}{3}M_{1}}-\frac{2M_{2}}{\pi}\frac{r_{0}(2-r_{0})}{(1-r_{0})^{2}}\\
&&-2M_{1}r_{0}(2+r_{0})\big[\log4(1+r_{0})- \log
r_{0}\big]\bigg\}\\
&=&\frac{2M_{2}}{\pi}\frac{r_{0}^{2}(2-r_{0})}{(1-r_{0})^{2}}.
\end{eqnarray*}
Hence $f(\mathbb{D}_{r_{0}})$ contains an univalent disk
$\mathbb{D}_{R_{0}}$ with
$$R_{0}\geq\frac{2M_{2}}{\pi}\frac{r_{0}^{2}(2-r_{0})}{(1-r_{0})^{2}}.$$
The proof of this theorem is complete. \qed

\subsection*{Proof of Corollary \ref{cor-1}}
For   $z_{1}, z_{2}\in\mathbb{D}_{r_{0}}$,  by (\ref{eq-bi-1}), we
see that there is a positive constant $L_{1}$ such that
$$L_{1}|z_{1}-z_{2}|\leq|f(z_{1})-f(z_{2})|,$$ where $r_{0}$ satisfies the
following equation
$$\frac{1}{\frac{4}{\pi}M_{2}+\frac{2}{3}M_{1}}-\frac{4M_{2}}{\pi}\frac{r_{0}(2-r_{0})}{(1-r_{0})^{2}}-2M_{1}\big[\log4(1+r_{0})- \log r_{0}\big](2+r_{0})r_{0}=0.$$

On the other hand, for   $z_{1}, z_{2}\in\mathbb{D}_{r_{0}}$, we use
Theorem \ref{thm2.0} to get
\begin{eqnarray*}
|f(z_{2})-f(z_{1})|&=&\left|\int_{[z_{1},z_{2}]}df(z)\right|\\&\leq&\int_{[z_{1},z_{2}]}\|D_{f}(z)\||dz|\\
&\leq&\int_{[z_{1},z_{2}]}\left(\frac{4M_{2}}{\pi}\frac{1}{1-r_{0}^{2}}+\frac{2}{3}M_{1}\right)|dz|\\
&=&\left(\frac{4M_{2}}{\pi}\frac{1}{1-r_{0}^{2}}+\frac{2}{3}M_{1}\right)|z_{1}-z_{2}|,
\end{eqnarray*}
where $[z_{1},z_{2}]$ is the segment from $z_{1}$ to $z_{2}$ with
the endpoints $z_{1}$ and $z_{2}$. Therefore,  $f$ is bi-Lipschitz
in $\mathbb{D}_{r_{0}}$.  \qed

\normalsize

\end{document}